\documentclass[10pt]{article}
\usepackage{amsfonts}
\usepackage{hyperref}
\usepackage{mathrsfs,enumerate,amssymb,amsmath,color,lineno}
\usepackage{indentfirst}
\usepackage{enumitem}
\usepackage{amsmath}
\usepackage{amssymb}
\usepackage[amsmath,thmmarks]{ntheorem}
\usepackage{graphicx}
\usepackage{tikz}

\newtheorem{definition}{\bf Definition}[section]
\newtheorem{lemma}{\bf Lemma}[section]
\newtheorem{theorem}{\bf Theorem}[section]
\newtheorem{remark}{\bf Remark}[section]
\newtheorem{corollary}{\bf Corollary}[section]
\newtheorem{example}{\bf Example}[section]

\newtheorem{proposition}{\bf Proposition}[section]

\begin{document}
\setcounter{page}{1}

\title{{\textbf{Further results on fuzzy negations and implications induced by fuzzy conjunctions and disjunctions}}\thanks {Supported by
the National Natural Science Foundation of China (No.12471440)}}
\author{Xin-Tong Zhang\footnote{\emph{E-mail address}: 1826124204@qq.com}, Xue-ping Wang\footnote{Corresponding author. xpwang1@hotmail.com; fax: +86-28-84761502},\\
\emph{School of Mathematical Sciences, Sichuan Normal University,}\\
\emph{Chengdu 610066, Sichuan, People's Republic of China}}

\newcommand{\pp}[2]{\frac{\partial #1}{\partial #2}}
\date{}
\maketitle
\begin{quote}
{\bf Abstract}  In this article, we deeply investigate some properties of fuzzy negations induced from fuzzy conjunctions (resp. disjunctions), which are then applied to characterizing the fuzzy negations. We further use the obtained characterization of fuzzy negations to explore some properties of $(D,N)$-implications generated from fuzzy disjunctions and negations. We finally describe $(D,N)$-implications (resp. continuous $(D,N)$-implications) generated from fuzzy disjunctions and negations. 

\end{quote}

\section{Introduction}
The indispensable role of fuzzy logic connectives is out of discussion in fuzzy logical systems and different logical connectives determine different logical systems. It is well-known that there are basic four types of logical connectives: a conjunction, a disjunction, a negation, and an implication. They are usually interpreted as corresponding operations in the theory of fuzzy logic systems. Specific speaking, a triangular norm (t-norm for short) is used to model the conjunction connective. Dually, a triangular conorm (t-conorm for short) is used to model the disjunction connective \cite{EP2000}. The negation connective is interpreted as a fuzzy negation and the implication connective is interpreted as a fuzzy implication \cite{MB2008}. Consequently, t-norms, t-conorms, fuzzy negations, fuzzy implications, and even more general aggregation operators such as uninorms, semi-uninorms and copulas \cite{CA2006} have formed several branches in the field of fuzzy logic systems. In addition, fuzzy logic connectives also have important applications in fuzzy system analyses, fuzzy integrals, fuzzy measures \cite{BB1993}, fuzzy decision-makings \cite{JF1994}, data mining, and probability metric spaces \cite{BS1983} etc.,
which lead to the construction of fuzzy negations and implications being an important content of fuzzy logic systems. For instance, Baczy\'{n}ski and Jayaram \cite{MB2008} and De Lima et al. \cite{AA2020} applied t-norms, t-conorms, uninorms, generators, ordinal sum, etc., to construct fuzzy negations and implications. Because the associativity of t-norms (resp. t-conorms) imposes severe limitations on their applicability as fuzzy conjunctions (resp. disjunctions), some authors like Durante et al. \cite{FD2007} consider to relax the conditions and use semi-copulas, quasi-copulas and copulas to induce and characterize residual implications. In particular, Jayaram \cite{BJ2017} investigated some properties of t-subnorms with strong associated negation and showed the conditions under which such t-subnorms are necessarily t-norms. After Fodor and Keresztfalvi \cite{JT1995} proposed fuzzy conjunctions and disjunctions which only satisfy the monotonicity and boundary conditions, Kr\'{o}l \cite{AK2011} explored basic properties of both residual implications generated from the fuzzy conjunctions and fuzzy conjunctions induced by fuzzy implications. Moreover, based on Baczy\'{n}ski and Jayaram \cite{MB2007, BJ2008}, Kr\'{o}l\cite{AK2012} also defined fuzzy negations induced by fuzzy conjunctions (resp. disjunctions) and obtained their basic properties. Kr\'{o}l \cite{AK20111} further defined $(D,N)$-implications generated by fuzzy disjunctions and negations and gave us some basic properties of $(D,N)$-implications. Later, Ma and Zhou \cite{QM2017} provided some characterizations of $(D,N)$-implications through strengthening the conditions of fuzzy disjunctions. To the best of our knowledge, there is currently no comprehensive characterization that applies to both fuzzy negations induced from fuzzy conjunctions (resp. disjunctions) and fuzzy implications induced by fuzzy disjunctions and negations in a unified theoretical framework. This article will consider the above problem.

The rest of this article is organized as follows. In Section 2, we present some necessary definitions and known results that will be used in the sequel. In Section 3, we deeply investigate some properties of fuzzy negations induced from fuzzy conjunctions (resp. disjunctions), which are then employed for characterizing the fuzzy negations. In Section 4, we first consider some properties of $(D,N)$-implications generated from fuzzy disjunctions and negations. Then we describe $(D,N)$-implications (resp. continuous $(D,N)$-implications) generated from fuzzy disjunctions and negations by imposing a slightly strengthened condition on fuzzy implications (resp. negations). A conclusion is drawn in Section 5.

\section{Preliminaries}
In this section, we recall some basic concepts and known results.
\begin{definition}[\cite{EP2000}]
\emph{A binary operator $T:[0,1]^2\rightarrow[0,1]$ is called a t-norm (resp. t-conorm), if it is non-decreasing in both variables, commutative, associative and has $1$ (resp. $0$) as the neutral element.}
\end{definition}

\begin{definition}[\cite{MB2008}]\label{def2.2}
\emph{ An operator $N:[0,1]\rightarrow [0, 1]$ is called a fuzzy negation if}

\emph{(N1) $N(0)=1,N(1)=0$,}

\emph{(N2) $N$ is non-increasing.}

\noindent \emph{Further, a fuzzy negation $N$ is strict if, in addition,}

\emph{(N3) $N$ is strictly decreasing,}

\emph{(N4) $N$ is continuous.}\\
\emph{Moreover, a fuzzy negation $N$ is called strong if it is an involution, i.e.,}

\emph{(N5) $N(N(x))=x$ for any $x\in[0,1]$.}
\end{definition}

\begin{lemma}[\cite{MB2008}]\label{lem2.1}
  Every strong fuzzy negation is strict.
\end{lemma}
\begin{lemma}[\cite{MB2008}]\label{lem2.2}
If $N_1$,$N_2$ are fuzzy negations such that $N_1\circ N_2=\emph{id}_{[0,1]}$, then $N_1$ is a continuous fuzzy negation and $N_2$ is a strictly decreasing fuzzy negation.
\end{lemma}
\begin{example}[\cite{MB2008, EP2000}]
\emph{The following are some examples of fuzzy negations:}
\begin{enumerate}[label=(\arabic*)]
	\item \emph{The least negation: $N_{G_1}(x)=\begin{cases}
                               1 & \mbox{if } x=0, \\
                               0 & \mbox{if } x\in(0,1].
                             \end{cases}$}

\item \emph{The greatest negation: $N_{G_2}(x)=\begin{cases}
                               1 & \mbox{if } x\in[0,1), \\
                               0 & \mbox{if } x=1.
                             \end{cases}$}

\item \emph{The standard negation: $N_S(x)=1-x$.}
\end{enumerate}
\end{example}

\begin{definition}[\cite{MB2008}]
\emph {A binary operator $I:[0, 1]^2\rightarrow [0, 1]$ is called a fuzzy implication if for all $x,x_1,x_2,y,y_1,y_2\in[0,1]$, it satisfies the following conditions:}

\emph{(I1) If $x_1\leq x_2$, then $I(x_1,y)\geq I(x_2,y)$, i.e., $I(\cdot,y)$ is non-increasing,}

\emph{(I2) If $y_1 \leq y_2$, then $I(x,y_1) \leq I(x,y_2)$, i.e., $I(x,\cdot)$ is non-decreasing,}

\emph{(I3) $I(0,0)=1$,}

\emph{(I4) $I(1,1)=1$,}

\emph{(I5) $I(1,0)=0$.}

\emph{The set of all fuzzy implications will be denoted by $\mathcal{FI}$. A fuzzy implication $I$ is said to satisfy}

\emph{(i) the neutral property if}
\begin{equation}\tag{NP}\label{NP}
  I(1,y)=y\emph{ for any }y\in[0,1];
\end{equation}

\emph{(ii) the exchange principle if}
\begin{equation}\tag{EP}\label{EP}
  I(x,I(y,z))=I(y,I(x,z))\emph{ for any }x,y,z\in[0,1];
\end{equation}

\emph{(iii) the identity principle if
\begin{equation}\tag{IP}\label{IP}
  I(x,x)=1\mbox{ for any }x\in[0,1];
\end{equation}}

\emph{(iv) the ordering property if
\begin{equation}\tag{OP}\label{OP}
  I(x,y)=1\Leftrightarrow x\leq y\mbox{ for any }x,y\in[0,1].
\end{equation}}
\end{definition}
\begin{lemma}[\cite{MB2008}]\label{lem2.3}
If a fuzzy implication $I$ fulfils (OP), then $I$ satisfies (IP).
\end{lemma}

\begin{definition}[\cite{JD2010}]
\emph{(i) A binary operator $C:[0,1]^2\rightarrow [0,1]$ is called a fuzzy conjunction if it is non-decreasing with respect to each variable,}
\begin{equation*}
  C(1,1)=1\emph{ and }C(0,0)=C(1,0)=C(0,1)=0.
\end{equation*}
\emph{(ii) A binary operator $D:[0,1]^2\rightarrow [0,1]$ is called a fuzzy disjunction if it is non-decreasing with respect to each variable,
\begin{equation*}
  D(0,0)=0\mbox{ and }D(1,1)=D(1,0)=D(0,1)=1.
\end{equation*}}
\end{definition}

\begin{lemma}[\cite{AK2011}]\label{lem2.4}
  A fuzzy conjunction (resp. disjunction) has an absorbing element $0$ (resp. $1$).
\end{lemma}
\begin{definition}[\cite{MB2008}]
\emph {Let $I$ be a fuzzy implication and $N$ be a fuzzy negation. Then $I$ is said to satisfy\\
(i) law of contraposition w.r.t. $N$, if
\begin{equation}\tag{CP($N$)}\label{CP($N$)}
	I(x,y)=I(N(y),N(x))\mbox{ for any }x,y\in[0,1];
\end{equation}
(ii) law of left contraposition w.r.t. $N$, if
\begin{equation}\tag{L-CP($N$)}\label{L-CP($N$)}
	I(N(x),y)=I(N(y),x)\mbox{ for any }x,y\in[0,1];
\end{equation}
(iii) law of right contraposition w.r.t. $N$, if
\begin{equation}\tag{R-CP($N$)}\label{R-CP($N$)}
	I(x,N(y))=I(y,N(x))\mbox{ for any }x,y\in[0,1].
\end{equation}}
\end{definition}
\section{Fuzzy negations induced from fuzzy conjunctions (resp. disjunctions)}
Based on the fuzzy negations induced from fuzzy conjunctions (resp. disjunctions) introduced by Kr\'{o}l \cite{AK2012}, in this section, we explore their properties deeply, which are then applied to characterizing the fuzzy negations.

We first cite the following definition from \cite{AK2012}.
\begin{definition}
\emph{(i) Let $C$ be a fuzzy conjunction. A function $N_C:[0,1]\rightarrow [0,1]$ defined as
\begin{equation}\tag{3.1}\label{3.1}
N_C(x)=\sup\{t\in[0,1]\mid C(x,t)=0\}\mbox{ for any } x\in[0,1]
\end{equation}
is called the natural negation of $C$ or the negation induced by $C$.}\\
\emph{(ii) Let $D$ be a fuzzy disjunction. A function $N_D:[0,1]\rightarrow [0,1]$ defined as
\begin{equation}\tag{3.2}\label{3.2}
N_D(x)=\inf\{t\in[0,1]\mid D(x,t)=1\}\mbox{ for any } x\in[0,1]
\end{equation}
is called the natural negation of $D$ or the negation induced by $D$.}
\end{definition}

\begin{example}[\cite{MB2008, AK2012}]\label{exa3.1}
\emph{Table 1 gives the natural negations of some fuzzy conjunctions, respectively. In Table 1, the symbol ``-" means that the function $N_C$ is not a fuzzy negation.}
\begin{table}[!h]
  \centering
  \caption{Fuzzy negations induced by fuzzy conjunctions}
  \begin{tabular}{|c|c|}
    \hline
    Fuzzy conjunction $C$ & Fuzzy negation $N_C$ \\
    \hline
    $C_0(x,y)=\begin{cases}
                1 & \mbox{if } (x, y)=(1, 1), \\
                0 & \mbox{otherwise}
              \end{cases}$ & $-$ \\
    \hline
    $C_1(x,y)=\begin{cases}
                0 & \mbox{if } x=0\mbox{ or }y=0, \\
                1 & \mbox{otherwise}
              \end{cases}$ & $N_{G_1}$ \\
    \hline
    $C_2(x,y)=\begin{cases}
                y & \mbox{if } x=1, \\
                0 & \mbox{otherwise}
              \end{cases}$ & $N_{G_2}$ \\
    \hline
    $C_3(x,y)=\begin{cases}
                x & \mbox{if } y=1, \\
                0 & \mbox{otherwise}
              \end{cases}$ & $-$ \\
    \hline
    $C_4(x,y)=\begin{cases}
                0 & \mbox{if } x+y\leq1, \\
                y & \mbox{otherwise}
              \end{cases}$ & $N_S$ \\
    \hline
    $C_5(x,y)=\begin{cases}
                0 & \mbox{if } y=0, \\
                x & \mbox{otherwise}
              \end{cases}$ & $N_{G_1}$ \\
    \hline
    $T_M(x,y)=\min\{x,y\}$ & $N_{G_1}$ \\
    \hline
    $T_P(x,y)=xy$ & $N_{G_1}$ \\
    \hline
    $T_L(x,y)=\max\{x+y-1,0\}$ & $N_S$ \\
    \hline
    $T_D(x,y)=\begin{cases}
                0 & \mbox{if } (x,y)\in[0,1)^2, \\
                \min\{x,y\} & \mbox{otherwise}
              \end{cases}$ & $N_{G_2}$ \\
    \hline
    $T_{nM}(x,y)=\begin{cases}
                 0 & \mbox{if } x+y\leq1, \\
                 \min\{x,y\} & \mbox{otherwise}
               \end{cases}$ & $N_S$ \\
    \hline
  \end{tabular}
  \end{table}
\end{example}

Then we have the following two propositions.
\begin{proposition}
If a fuzzy conjunction $C$ is left-continuous, then

\noindent(i) for some $x,y\in[0,1]$,
\begin{equation}\tag{3.3}\label{3.3}
C(x,y)=0\Longleftrightarrow N_C(x)\geq y
\end{equation}

\noindent(ii)
$N_C(x)=\max\{t\in[0,1]\mid C(x,t)=0\}\mbox{ for any } x\in[0,1]$.

\noindent(iii) $N_C$ is left-continuous.
\end{proposition}

\noindent\begin{proof}
(i) Assume that $C$ is a left-continuous fuzzy conjunction. If $C(x,y)=0$ for some $x,y\in[0,1]$, then $y\in\{t\in[0,1]\mid C(x,t)=0\}$ and hence $N_C(x)\geq y$ by Eq. \eqref{3.1}.

 Now, suppose that $N_C(x)\geq y$ for some $x,y\in[0,1]$. There are two cases as follows:

\noindent Case 1. If $N_C(x)>y$, then there exists some $t>y$ such that $C(x,t)=0$, which implies that $C(x,y)\leq C(x,t)=0$ by the monotonicity of $C$, i.e., $C(x,y)=0$.

\noindent Case 2. If $N_C(x)=y$, then either $y\in\{t\in[0,1]\mid C(x,t)=0\}$ or $y\notin\{t\in[0,1]\mid C(x,t)=0\}$. In the first case, $C(x,y)=0$. In the other case, there exists an increasing sequence $(t_i)_{i\in\mathbb{N}}$ such that $C(x,t_i)=0$ for all $i\in\mathbb{N}$ and $\lim_{i\to\infty}t_i=y$. By the left-continuous of $C$ we get
\begin{equation*}
 C(x,y)=C(x,\lim_{i\to\infty}t_i)=\lim_{i\to\infty}C(x,t_i)=\lim_{i\to\infty}0=0,
\end{equation*}
\noindent a contradiction.

 From Cases 1 and 2, we have that $C(x,y)=0$ for some $x,y\in[0,1]$.

\noindent(ii) Because $N_C(x)\leq N_C(x)$ for all $x\in[0,1]$,  from (i) we have $C(x,N_C(x))=0$, which means $N_C(x)\in\{t\in[0,1]\mid C(x,t)=0\}$. Therefore, $N_C(x)=\max\{t\in[0,1]\mid C(x,t)=0\}\mbox{ for any } x\in[0,1]$.

\noindent(iii) Assume that $N_C$ is not left-continuous at $x_0\in(0,1]$. Then there exist $a,b\in[0,1]$ with $a>b$ such that $N_C(x)\geq a$ for any $x<x_0$ and $N_C(x_0)=b$ since $N_C$ is non-increasing. Thus by Eq.\eqref{3.3} we get $C(x,a)=0$ for any $x<x_0$. Since $C$ is left-continuous, we have $0=\lim_{x\to x_0}C(x,a)=C(x_0,a)$. Again from Eq.\eqref{3.3} we obtain $b=N_C(x_0)\geq a$, a contradiction.
\end{proof}

\begin{remark}\label{rem3.1}
\emph{From the proof of Proposition 3.1 (i), without the left-continuity of $C$, the following statements hold:}

\noindent \emph{(i) If $C(x,y)=0$ for some $x,y\in[0,1]$, then $y\leq N_C(x)$.}

\noindent\emph{(ii) If $y<N_C(x)$ for some $x,y\in[0,1]$, then $C(x,y)=0$.}\\
\emph{Thus:}

\noindent\emph{(i$'$) If $y>N_C(x)$ for some $x,y\in[0,1]$, then $C(x,y)>0$.}

\noindent\emph{(ii$'$) If $C(x,y)>0$ for some $x,y\in[0,1]$, then $y\geq N_C(x)$.}
\end{remark}

\begin{proposition}\label{prop3.2}
\emph{Let $C$ be a commutative fuzzy conjunction and $N_C$ be a fuzzy negation. Then the following two statements hold:}

\noindent(i) If $N_C$ is continuous, then it is strong.

\noindent(ii) If $N_C$ is not continuous, then it is not strictly decreasing.
\end{proposition}

\noindent\begin{proof}
\noindent(i) We first show that $N_C$ is strictly monotone since $N_C$ is continuous. Assume that $N_C$ is not strictly monotone. Then there exist $x,y\in[0,1]$ with $0<x<y<1$ such that
\begin{equation*}
  N_C(x)=N_C(y)=p.
\end{equation*}

If $p=0$, then $N_C(x)=0$. This implies that $C(x,\varepsilon)>0$ for an arbitrary $\varepsilon\in (0,1]$. By the commutativity of $C$, $C(\varepsilon,x)=C(x,\varepsilon)>0$. Thus by Remark \ref{rem3.1} (ii$'$) we have $N_C(\varepsilon)\leq x$. On the other hand, $\lim_{\varepsilon\to 0}N_C(\varepsilon)=N_C(\lim_{\varepsilon\to 0}\varepsilon)=N_C(0)=1$ since $N_C$ is continuous. Therefore, $1\leq x$, a contradiction.

If $p=1$, then $N_C(y)=1$. We claim that $C(y,1-\varepsilon)=0$ for an arbitrary $\varepsilon\in [0,1)$. Otherwise, $C(y,1-\varepsilon)>0$ for some $\varepsilon\in [0,1)$. Then by Remark \ref{rem3.1} (ii$'$) we have $1-\varepsilon\geq N_C(y)=1$, a contradiction. Since $C$ is commutative, $C(1-\varepsilon,y)=C(y,1-\varepsilon)=0$. Thus by Remark \ref{rem3.1} (i) we have $N_C(1-\varepsilon)\geq y$. On the other hand, $\lim_{\varepsilon\to 0}N_C(1-\varepsilon)=N_C(\lim_{\varepsilon\to 0}(1-\varepsilon))=N_C(1)=0$ since $N_C$ is continuous. Therefore, $0=N_C(1)\geq y$, a contradiction.

 Now assume $p\in(0,1)$. Since $N_C(z)=p$ for any $z\in(x,y)$, analogously we have
\begin{equation*}
  C(z,p-\varepsilon)=C(p-\varepsilon,z)=0
\end{equation*}
for any arbitrary $p-\varepsilon\in (0,1)$ and
\begin{equation*}
  C(z,p+\varepsilon)=C(p+\varepsilon,z)>0
\end{equation*}
for any arbitrary $p+\varepsilon\in (0,1)$. Once again, by Remark \ref{rem3.1} (i) and (ii$'$) we get
\begin{equation*}
  N_C(p+\varepsilon)\leq z\leq N_C(p-\varepsilon).
\end{equation*}
Thus $\lim_{\varepsilon\to 0}N_C(p+\varepsilon)\leq z\leq \lim_{\varepsilon\to 0}N_C(p-\varepsilon)$. This follows  $N_C(p)=z$ since $N_C$ is continuous, i.e., $N_C(p)=z$ for every $z\in(x,y)$, which contradicts the fact that $N_C$ is a function.

Consequently, $N_C$ is strictly monotone. Therefore, $N_C$ is strict.

To complete the proof, we need to show $N_C(N_C(x))=x$ for all $x\in[0,1]$. Since $N_C$ is strict, it suffices to prove $N_C(N_C(N_C(x)))=N_C(N_C(x))$ for all $x\in[0,1]$. First, it is clear that $N_C(N_C(N_C(x)))=N_C(N_C(x))$ for any $x\in\{0,1\}$. Secondly, for any $x\in(0,1)$ and an arbitrary small $\varepsilon>0$ we have $x-\varepsilon<x<x+\varepsilon$. This follows that
\begin{equation}\tag{3.4}\label{3.4}N_C(x+\varepsilon)<N_C(x)<N_C(x-\varepsilon),\end{equation} thus \begin{equation}\tag{3.5}\label{3.5}N_C(N_C(x-\varepsilon))<N_C(N_C(x))<N_C(N_C(x+\varepsilon))\end{equation} since $N_C$ is strict. Then by Eq. \eqref{3.4} and Remark \ref{rem3.1} (ii) we obtain
\begin{equation*}
C(x,N_C(x+\varepsilon))=C(N_C(x+\varepsilon),x)=0.
\end{equation*}
Thus by Remark \ref{rem3.1} (i) we have $x\leq N_C(N_C(x+\varepsilon))$, which implies that
\begin{equation*}
 N_C(x)\geq N_C(N_C(N_C(x+\varepsilon)))
\end{equation*}
since $N_C$ is strictly decreasing. Hence, by the continuity of $N_C$ we get
\begin{equation}\tag{3.6}\label{3.6}
 N_C(x)\geq \lim_{\varepsilon\to 0}N_C(N_C(N_C(x+\varepsilon)))= N_C(N_C(N_C(x))).
\end{equation}
On the other hand, by Eq. \eqref{3.5} and Remark \ref{rem3.1} (ii), we have
\begin{equation*}
 C(N_C(N_C(x-\varepsilon)),N_C(x))=C(N_C(x),N_C(N_C(x-\varepsilon)))=0.
\end{equation*}
Thus by Remark \ref{rem3.1} (i) we have $N_C(x)\leq N_C(N_C(N_C(x-\varepsilon)))$. Hence
\begin{equation*}
 N_C(x)\leq \lim_{\varepsilon\to 0}N_C(N_C(N_C(x-\varepsilon)))= N_C(N_C(N_C(x)))
\end{equation*}since $N_C$ is continuous, which together with Eq. \eqref{3.6} implies
\begin{equation*}
 N_C(x)=N_C(N_C(N_C(x))).
\end{equation*}
Therefore, $N_C(N_C(N_C(x)))=N_C(N_C(x))$ for all $x\in[0,1]$.

\noindent(ii) If $N_C$ is not continuous at $p\in[0,1]$, then by the monotonicity of $N_C$, there exist $x,y\in[0,1]$ with $x< y$ such that
\begin{equation*}
\begin{aligned}
& x=\begin{cases}
\lim_{t\to p^+}N_C(t) & \hbox{if }\ p<1,\\
0 & \hbox{if }\ p=1\\
\end{cases}\mbox{ and }
&&
& y=\begin{cases}
\lim_{t\to p^-}N_C(t) & \hbox{if }\ p>0,\\
1 & \hbox{if }\ p=0.\\
\end{cases}
\end{aligned}
\end{equation*}
Now, we consider the following two cases:

\noindent Case 1.  If $N_C(p)=x$, then $p>0$. It is obvious that $N_C(p)<z$ for any $z\in(x,y)$. Thus $C(p,z)>0$ by Remark \ref{rem3.1} (i$'$), which follows that $C(z,p)=C(p,z)>0$ since $C$ is commutative. By using Remark \ref{rem3.1} (ii$'$), we get $N_C(z)\leq p$. On the other hand, it is clearly that $N_C(p-\varepsilon)\geq y$ for an arbitrary small $\varepsilon>0$. Thus $z<y\leq N_C(p-\varepsilon)$. It follows from Remark \ref{rem3.1} (ii) and the commutativity of $C$ that $C(z,p-\varepsilon)=C(p-\varepsilon,z)=0$. Then from Remark \ref{rem3.1} (i) we have $N_C(z)\geq p-\varepsilon$. Thus $N_C(z)\geq \lim_{\varepsilon\to 0}(p-\varepsilon)=p$.
Therefore, $N_C(z)=p$ for every $z\in(x,y)$, i.e., $N_C$ is not strictly decreasing.

\noindent Case 2.  If $N_C(p)=z'$ with $x<z'\leq y$, then $p<1$. We now claim that $N_C$ is constant on the interval $(x,z')$. Indeed, let $z\in(x,z')$. Then $z<z'=N_C(p)$. Thus by Remark \ref{rem3.1} (ii) and the commutativity of $C$, we have $C(p,z)=C(z,p)=0$. Hence $N_C(z)\geq p$ by Remark \ref{rem3.1} (i). If $N_C(z)>p$, then $N_C(z)>p+\varepsilon$ for some $\varepsilon>0$. It follows from Remark \ref{rem3.1} (ii) that $C(p+\varepsilon,z)=C(z,p+\varepsilon)=0$. Again by Remark \ref{rem3.1} (i), we get $N_C(p+\varepsilon)\geq z$. Since $N_C$ is non-increasing, we obtain
\begin{equation*}
  x\geq N_C(p+\varepsilon)\geq z,
\end{equation*}
\noindent a contradiction. Therefore, $N_C(z)=p$ for every $z\in(x,z')$, i.e., $N_C$ is not strictly decreasing.
\end{proof}

The following example illustrates that in Proposition \ref{prop3.2} the commutativity of $C$ is necessary.
\begin{example}
\emph{\noindent (i) Consider the fuzzy conjunction $C:[0,1]^2\rightarrow [0,1]$ defined by
\begin{equation*}
C(x,y)=\max\{x+\sqrt{y}-1,0\}.
\end{equation*}
Thus $N_C(x)=(1-x)^2$. Clearly, $C$ is not commutative and $N_{C}$ is continuous but not strong.}

\emph{\noindent (ii) Consider the fuzzy conjunction $C:[0,1]^2\rightarrow [0,1]$ defined by}
\begin{equation*}
 C(x,y)=\begin{cases}
          \max\{0.5x+y-1,0\} & \emph{if } x\in[0,0.5), \\
          \max\{x+y-1,0\} & \emph{otherwise}.
        \end{cases}
\end{equation*}

\emph{Thus} \begin{equation*}
      N_{C}(x)=\begin{cases}
               1-0.5x & \emph{if } x\in[0,0.5), \\
               1-x & \emph{otherwise}.
             \end{cases}
    \end{equation*}

\emph{Clearly, $C$ is not commutative and $N_{C}$ is not continuous but strictly decreasing.}

\emph{\noindent(iii) From Example \ref{exa3.1}, we know that $N_{T_{nM}}(x)=1-x$ is a strong fuzzy negation. Analogously, we can construct a series of fuzzy conjunctions $C$ as follows:}

\begin{equation*}
C(x,y)=\begin{cases}
         0 & \emph{if } x+y\leq 1, \\
         \ast & \emph{otherwise},
       \end{cases}
\end{equation*}
\emph{in which the $\ast$ only needs to satisfy that $C(1,1)=1$, $C(x,y)>0$ when $x+y> 1$ and $C$ is non-decreasing with respect to each variable. Therefore, $N_C(x)=1-x$ and $C$ need not be commutative, for instance, the fuzzy conjunction $C_4$ of Example \ref{exa3.1}.}
\end{example}

The following theorem gives a characterization of fuzzy negations induced by fuzzy conjunctions.
\begin{theorem}
Let $C$ be a commutative fuzzy conjunction and $N_C$ be a fuzzy negation. Then the following statements are equivalent:

\noindent(i) $N_C$ is strictly decreasing.

\noindent(ii) $N_C$ is continuous.

\noindent(iii) $N_C$ is strict.

\noindent(iv) $N_C$ is strong.
\end{theorem}

\noindent\begin{proof}
By Definition \ref{def2.2}, it is obvious that (iv) implies  (iii) and (iii) deduces (i). From Proposition \ref{prop3.2} (i) we infer (iv) from (ii). Finally, by Proposition \ref{prop3.2} (ii) we know that (i) implies (ii).
\end{proof}

 Dually, for a fuzzy disjunction $D$ we have the following results.
\begin{proposition}\label{pro3.3}
If a fuzzy disjunction $D$ is right-continuous, then\\
(i) for some $x,y\in[0,1]$,
\begin{equation}\tag{3.7}\label{3.7}
D(x,y)=1\Longleftrightarrow N_D(x)\leq y
\end{equation}
(ii) $N_D(x)=\inf\{t\in[0,1]\mid D(x,t)=1\}\mbox{ for any } x\in[0,1],$

\noindent$(iii)$ $N_D$ is right-continuous.
\end{proposition}
\begin{remark}\label{rem3.2}
	\emph{Dual to Remark \ref{rem3.1}, for a fuzzy disjunction $D$ without the right-continuity, the following statements hold:}
	
	\noindent\emph{(i) If $D(x,y)=1$ for some $x,y\in[0,1]$, then $y\geq N_D(x)$.}
	
	\noindent\emph{(ii) If $y>N_D(x)$ for some $x,y\in[0,1]$, then $D(x,y)=1$.}
\end{remark}
\begin{proposition}\label{the3.2}
Let $D$ be a commutative fuzzy disjunction and $N_D$ be a fuzzy negation. Then the following two statements hold:

\noindent(i) If $N_D$ is continuous, then it is strong.

\noindent(ii) If $N_D$ is not continuous, then it is not strictly decreasing.
\end{proposition}

\begin{theorem}\label{thm3.2}
Let $D$ be a commutative fuzzy disjunction and $N_D$ be a fuzzy negation. Then the following statements are equivalent:

\noindent(i) $N_D$ is strictly decreasing.

\noindent(ii) $N_D$ is continuous.

\noindent(iii) $N_D$ is strict.

\noindent(iv) $N_D$ is strong.
\end{theorem}
\section{Fuzzy implications induced by fuzzy disjunctions and negations}

In classical logic, the implication operation can be expressed using negation and disjunction, i.e.,
\begin{equation*}
	p\rightarrow q\equiv \neg p\vee q.
\end{equation*}
Inspired by this, replacing $\neg$ with a fuzzy negation $N$, and replacing $\vee$ with a t-conorm $S$, a uninorm $U$, or a grouping function $G$, we can obtain various classes of fuzzy implications: $(S,N)$-implications (\cite{ MB2007,MB2008}), $(U,N)$-implications (\cite{MB2009}) and $(G,N)$-implications (\cite{GD2014}). Up to now, these three classes of implications have been relatively well studied. In particular, because $S$, $U$, $G$ are all specific fuzzy disjunctions, Kr\'{o}l \cite{AK20111} incorporated them into a unified framework and introduced the concept of $(D,N)$-implications. Ma and Zhou \cite{QM2017} further investigate basic properties and characterizations of $(D,N)$-implications through strengthening the conditions of fuzzy disjunctions. In this section, we continue to explore some properties and characterizations of $(D,N)$-implications (resp. continuous $(D,N)$-implications) by imposing a slightly strengthened condition on fuzzy implications (resp. negations).

First, we need the following definition.
\begin{definition}[\cite{AK2012}]\label{def4.1}
\emph{(i) Let $D$ be a fuzzy disjunction and $N$ a fuzzy negation. The pair $(D, N)$ is said to satisfy the law of excluded middle (LEM for short) if both
\begin{equation}\tag{LEM1}\label{LEM1}
 D(N(x),x)=1\mbox{ for any }x\in[0,1]
\end{equation}
and
\begin{equation}\tag{LEM2}\label{LEM2}
 D(x,N(x))=1\mbox{ for any } x\in[0,1].
\end{equation}
(ii) Let $C$ be a fuzzy conjunction and $N$ a fuzzy negation. The pair $(C, N)$ is said to satisfy the law of contradiction (LC for short) if both
\begin{equation}\tag{LC1}\label{LC1}
 C(N(x),x)=0\mbox{ for any } x\in[0,1]
\end{equation}
and
\begin{equation}\tag{LC2}\label{LC2}
 C(x,N(x))=0\mbox{ for any } x\in[0,1].
\end{equation}}
\end{definition}

Then we have the following lemma.
\begin{lemma}\label{lem4.1}
Let $D$ be a fuzzy disjunction and $N$ be a fuzzy negation. If $(D, N)$ satisfies (LEM), then

\noindent(i) $N_D(N(x))\leq x\mbox{ for any }x\in[0,1]$

\noindent(ii) $N(x)\geq N_D(x)\mbox{ for any }x\in[0,1]$
\end{lemma}
\begin{proof}
Assume that $(D, N)$ satisfies (LEM). It is obvious that by Remark \ref{rem3.2} (i) and (LEM1), $N_D(N(x))\leq x$ for any $x\in[0,1]$. Once again we have $N(x)\geq N_D(x)$ for any $x\in[0,1]$ by Remark \ref{rem3.2} (i) and (LEM2).
\end{proof}

\begin{remark}
\emph{The inverse of Lemma \ref{lem4.1} is not true in general. For example, let}
\begin{equation*}
D(x,y)=\begin{cases}
           1 & \emph{if } x+y>1, \\
           x^2+y^2 & \emph{otherwise}.
         \end{cases}
\end{equation*}
\emph{Then $N_D(x)=N_S(x)=1-x$ and $N_D(N_S(x))=x$ for all $x\in[0,1]$. However, the pair $(D,N_S)$ does not satisfy (LEM). Indeed, let $x=0.5$. Then}
$D(N_S(0.5),0.5)=D(0.5,N_S(0.5))=D(0.5,0.5)=0.5\neq 1$.
\end{remark}

Fortunately, we have the following proposition.
\begin{proposition}
Let $D$ be a right-continuous fuzzy disjunction and $N$ be a fuzzy negation. Then the following statements are equivalent:\\
(i) The pair $(D,N)$ satisfies (LEM).\\
(ii)  $N_D(N(x))\leq x$ and $N(x)\geq N_D(x)$ for any $x\in[0,1]$.
\end{proposition}

\noindent\begin{proof}
By Lemma \ref{lem4.1} it is enough to show that (ii) implies (i). Assume that $N_D(N(x))\leq x$ and $N(x)\geq N_D(x)$ for any $x\in[0,1]$. It is obvious that $D(N(x),x)=1$ and $D(x,N(x))=1$ by Eq.\eqref{3.7}, i.e., the pair $(D,N)$ satisfies (LEM).
\end{proof}

From Definition \ref{def4.1}, one easily see that (LEM) and (LC) are dual to each other. Then any result for (LEM) can dually be translated to (LC). Therefore, we have the following two results for (LC).
\begin{lemma}
Let $C$ be a fuzzy conjunction and $N$ be a fuzzy negation. If $(C, N)$ satisfies (LC), then\\
(i) $N_C(N(x))\geq x \mbox{ for any } x\in[0,1]$.\\
(ii) $N(x)\leq N_C(x) \mbox{ for any } x\in[0,1]$.
\end{lemma}
\begin{proposition}
Let $C$ be a left-continuous fuzzy conjunction and $N$ be a fuzzy negation. The following statements are equivalent:\\
(i) The pair $(C, N)$ satisfies (LC).\\
(ii) $N_C(N(x))\geq x$ and $N(x)\leq N_C(x)$ for any $x\in[0,1]$.
\end{proposition}

Before giving some equivalent descriptions of $I_{D,N}$ that satisfy (IP) (resp. (OP)), we need the following definition and lemma.
\begin{definition}[\cite{AK20111, QM2017}]\label{def4.2}
\emph{A binary operator $I:[0, 1]^2\rightarrow [0, 1]$ is called a  $(D,N)$-implication, denoted by $I_{D,N}$, if there exist a fuzzy disjunction $D$ and a fuzzy negation $N$ such that
\begin{equation}\tag{4.1}\label{4.1}
I(x,y)=D(N(x),y)\mbox{ for any }x,y\in[0,1].
\end{equation}}
\end{definition}

\begin{lemma}[\cite{AK20111, QM2017}]\label{lem4.3}
Let $D$ be a fuzzy disjunction and $N$ be a fuzzy negation. Then:\\
(i) $I_{D,N}\in\mathcal{FI}$;\\
(ii) $I_{D,N}$ satisfies (NP) if and only if $D(0,y)=y$ for any $y\in[0,1]$;\\
(iii) If $D$ is commutative and associative, then $I_{D,N}$ satisfies (EP).
\end{lemma}

The following proposition provides necessary and sufficient conditions for $I_{D,N}$ satisfying (IP) and (OP).
\begin{proposition}\label{pro4.4}
Let $D$ be a fuzzy disjunction, $N$ and $N_D$ be two fuzzy negations. Then:\\
(i) $I_{D,N}$ satisfies (IP) if and only if the pair $(D,N)$ satisfies (LEM1);\\
(ii) If $D$ is commutative, then $I_{D,N}$ satisfies (OP) if and only if the pair $(D,N)$ satisfies (LEM1), $N=N_D$ and both of them are strong fuzzy negations.
\end{proposition}
\begin{proof}
(i) $I_{D,N}$ satisfies (IP) $\Leftrightarrow$ $I_{D,N}(x,x)=1$ for any $x\in[0,1]$ $\Leftrightarrow$ $D(N(x),x)=1$ for any $x\in[0,1]$ $\Leftrightarrow$ $(D,N)$ satisfies (LEM1).

(ii) Suppose that $I_{D,N}$ satisfies (OP). Then $I_{D,N}$ satisfies (IP) by Lemma \ref{lem2.3}. Then by (i) we get that the pair $(D,N)$ satisfies (LEM1). Thus $N_D(N(x))\leq x$ for any $x\in[0,1]$ by Lemma \ref{lem4.1} (i). Assume that there is an $x_0$ such that $N_D(N(x_0))<x_0$. Then there exists a $y\in[0,1]$ such that $N_D(N(x_0))<y<x_0$. It follows from Remark \ref{rem3.2} (ii) that $I_{D,N}(x_0,y)=D(N(x_0),y)=1$. Thus $x_0\leq y$ since $I_{D,N}$ satisfies (OP), a contradiction. Hence $N_D(N(x))=x$ for any $x\in[0,1]$, which implies that $N_D$ is continuous and $N$ is strictly decreasing by Lemma \ref{lem2.1}. Because $D$ is commutative, by Theorem \ref{thm3.2} we get that $N_D$ is strong and $N=N_D$.

 Conversely, assume that the pair $(D,N)$ satisfies (LEM1), $N=N_D$ and both of them are strong fuzzy negations. Let $x,y\in[0,1]$ with $x\leq y$. Then by the monotonicity of $D$ we have $I_{D,N}(x,y)=D(N(x),y)\geq D(N(x),x)=1$, i.e., $I_{D,N}(x,y)=1$. Now, suppose that $I_{D,N}(x,y)=1$ for all $x,y\in[0,1]$, i.e., $D(N(x),y)=1$. Then by Remark \ref{rem3.2} (i) we obtain $y\geq N_D(N(x))$. Thus $y\geq x$ since $N=N_D$ is a strong fuzzy negation. Therefore, $I_{D,N}$ satisfies (OP).
\end{proof}

In particular, we have the following example.
\begin{example}
\emph{Let $D$ be a fuzzy disjunction. In Table 2 we provide some examples of $(D,N_D)$-implications.}
\begin{table}[!h]
  \centering
  \caption{$(D,N_D)$-implications induced by fuzzy disjunctions}\label{}
\begin{tabular}{|c|c|c|}
  \hline
   Fuzzy disjunction $D$ & Fuzzy negation $N_D$ & Implication $I_{D,N_D}$ \\
  \hline
  $D_1(x,y)=\begin{cases}
           y & \mbox{if } x=0, \\
           1 & \mbox{otherwise}
         \end{cases}$ & $N_{G_1}$ & $I_D(x,y)=\begin{cases}
                                                  1 & \mbox{if } x=0, \\
                                                  y & \mbox{otherwise}
                                                \end{cases}$ \\
  \hline
  $D_2(x,y)=\begin{cases}
         1 & \mbox{if } x+y\geq 1, \\
         y & \mbox{if } x+y<1
       \end{cases}$ & $N_S$ & $I_{GD}$ \\
  \hline
  $D_3(x,y)=\begin{cases}
         1 & \mbox{if } x+y\geq 1, \\
         x & \mbox{if } x+y<1
       \end{cases}$ & $N_S$ & $I_3(x,y)=\begin{cases}
                                          1 & \mbox{if } x\leq y, \\
                                          1-x & \mbox{if } x>y
                                        \end{cases}$ \\
  \hline
  $D_4(x,y)=\min\{1,x^2+y^2\}$& $N_4(x)=\sqrt{1-x^2}$ & $I_4(x,y)=\min\{1,1-x^2+y^2\}$ \\
  \hline
  $D_5(x,y)=\min\{1,x+\sqrt{y}\}$ & $N_5(x)=(1-x)^2$ & $I_5(x,y)=\min\{1,(1-x)^2+\sqrt{y}\}$ \\
  \hline
  $D_6(x,y)=\begin{cases}
              0 & \mbox{if } (x, y)=(0, 0), \\
              1 & \mbox{if } x=1\mbox{ or }y=1, \\
              x & \mbox{otherwise}
            \end{cases}$ & $N_{G_2}$ & $I_6(x,y)=\begin{cases}
                                                   0 & \mbox{if } x=1, y\in[0,1)\\
                                                   1 & \mbox{otherwise}
                                                 \end{cases}$ \\
  \hline
  $D_7(x,y)=\begin{cases}
              0 & \mbox{if }(x, y)=(0, 0), \\
              1 & \mbox{if } x=1 \mbox{ or }y=1, \\
              y & \mbox{otherwise}
            \end{cases}$ & $N_{G_2}$ & $I_{WB}$ \\
  \hline
\end{tabular}
\end{table}
\end{example}

Analogously, from Lemma \ref{lem4.3}, Proposition \ref{pro4.4} and Theorem \ref{thm3.2}, we have the following proposition.
\begin{proposition}\label{propo4.5}
Let $D$ be a fuzzy disjunction and $N_D$ be a fuzzy negation. Then:\\
(i) $I_{D,N_D}$ satisfies (NP) if and only if $D$ has a left neutral element $0$;\\
(ii) If $D$ is commutative and associative, then $I_{D,N_D}$ satisfies (EP);\\
(iii) $I_{D,N_D}$ satisfies (IP) if and only if the pair $(D,N_D)$ satisfies (LEM1);\\
(iv) $I_{D,N_D}$ satisfies (OP) if and only if the pair $(D,N_D)$ satisfies (LEM1) and $N_D$ is a strong fuzzy negation;\\
(v) If $D$ is commutative, then $I_{D,N_D}$ satisfies R-CP($N_D$). In addition, if $N_D$ is a continuous fuzzy negation, then $I_{D,N_D}$ satisfies L-CP($N_D$) and CP($N_D$).
\end{proposition}

In what follows, we consider a characterization of a (continuous) $(D,N)$-implication. We first need the following definition and three lemmas.
\begin{definition}[\cite{EP2000}]
\emph{Let $[a,b]$ and $[c,d]$ be two closed subintervals of the extended real line $[-\infty,+\infty]$, and $f:[a,b]\rightarrow[c,d]$ be a monotonic mapping. Define the pseudo-inverse $f^{(-1)}:[c,d]\rightarrow[a,b]$ as follows:
\begin{equation*}
  f^{(-1)}(y)=\sup\{x\in[a,b]\mid (f(x)-y)(f(b)-f(a))<0\}.
\end{equation*}}
\end{definition}

Note that if $f:[a,b]\rightarrow[c,d]$ is a non-increasing non-constant function, then for all $y\in[c,d]$,
\begin{equation*}
	f^{(-1)}(y)=\sup\{x\in[a,b]\mid f(x)>y\}.
\end{equation*}

\begin{lemma}[\cite{MB2008}]\label{lem5.2}
If $N$ is a continuous fuzzy negation, then the function $\aleph:[0,1]\rightarrow[0,1]$ defined by
\begin{equation}\tag{4.2}\label{4.2}
  \aleph(x)=\begin{cases}
              N^{(-1)}(x), & \mbox{if } x\in(0,1], \\
              1, & \mbox{if } x=0
            \end{cases}
\end{equation}
is a strictly decreasing fuzzy negation. Moreover,
\begin{equation*}
  \aleph^{(-1)}=N,
\end{equation*}
\begin{equation}\tag{4.3}\label{4.3}
  N\circ\aleph=\emph{id}_{[0,1]},
\end{equation}
\begin{equation}\tag{4.4}\label{4.4}
  \aleph\circ N\mid_{\emph{Ran}(\aleph)}=\emph{id}_{\emph{Ran}(\aleph)}.
\end{equation}
\end{lemma}

\begin{remark}[\cite{MB2008}]\label{rem4.2}
\emph{In Lemma \ref{lem5.2}, the pseudo-inverse $N^{(-1)}$ need not be a fuzzy negation even if $N$ is a continuous fuzzy negation. For instance, let $N: [0,1]\rightarrow[0,1]$ be defined by}
$$N(x)=\begin{cases}
        -2x+1 & \emph{if } x\in[0,0.5], \\
        0 & \emph{otherwise}.
      \end{cases}$$
\emph{Then $N^{(-1)}(x)=-0.5x+0.5$. Clearly, $N$ is continuous and $N^{(-1)}$ is not a fuzzy negation since $N^{(-1)}(0)=0.5\neq 1$.}
\end{remark}
\begin{lemma}[\cite{MB2008}]
	If a function $I:[0,1]^2\rightarrow[0,1]$ satisfies (I1), (I3) and (I5), then the function $N_I:[0,1]\rightarrow[0,1]$ defined by
	\begin{equation}\tag{4.5}\label{4.5}
		N_I(x)=I(x,0)\mbox{ for any } x\in[0,1]
	\end{equation}
	is a fuzzy negation.
\end{lemma}
\begin{lemma}[\cite{QM2017}]
Let $I$ be a fuzzy implication and $N_I$ be a continuous fuzzy negation. Define a function $D_I:[0,1]^2\rightarrow[0,1]$ by
\begin{equation}\tag{4.6}\label{4.6}
  D_I(x,y)=I(\aleph_I(x),y),
\end{equation}
then $D_I$ is a fuzzy disjunction.
\end{lemma}
\begin{remark}[\cite{QM2017}]\label{rem4.3}
\emph{(i) If $I$ is a fuzzy implication and $N_I$ is a strict fuzzy negation, then $\aleph_I=N^{-1}_I$. In this case,
\begin{equation*}
  D_I(x,y)=I(N^{-1}_I(x),y)\mbox{ for any }\ x,y\in[0,1].
\end{equation*}
(ii) If $I$ is a fuzzy implication and $N_I$ is a strong fuzzy negation, then $\aleph_I=N_I$. In this case,
\begin{equation*}
  D_I(x,y)=I(N_I(x),y)\mbox{ for any }x,y\in[0,1].
\end{equation*}}
\end{remark}
\begin{theorem}\label{the4.1}
For a function $I:[0,1]^2\rightarrow[0,1]$ the following statements are equivalent:\\
(i) $I$ is a $(D,N)$-implication in which $N$ is a continuous fuzzy negation and $D$ is a fuzzy disjunction with a right neutral element $0$.\\
(ii) $I$ is a fuzzy implication, $N_I$ is a continuous fuzzy negation and
\begin{equation}\tag{4.7}\label{4.7}
 I(x_1,y)=I(x_2,y) \mbox{ for any } y\in[0,1] \mbox{ whenever } \ N_I(x_1)=N_I(x_2)\mbox{ for some } x_1, x_2\in[0,1].
\end{equation}
Moreover, the representation Eq.\eqref{4.1} of $(D,N)$-implication is unique in this case.
\end{theorem}
\begin{proof}
(i)$\Rightarrow$(ii) Assume that $I$ is a $(D,N)$-implication in which $N$ is a continuous fuzzy negation and $D$ is a fuzzy disjunction with a right neutral element $0$. Then by Lemma \ref{lem4.3}, $I$ is a fuzzy implication. From Eqs. \eqref{4.1} and \eqref{4.5} we have $N_I(x)=I(x,0)=D(N(x),0)=N(x)$ for any $x\in [0,1]$ since $D$ is a fuzzy disjunction with a right neutral element $0$. Thus $N_I$ is a continuous fuzzy negation since $N$ is a continuous fuzzy negation. Now, suppose that $N_I(x_1)=N_I(x_2)$ for some $x_1,x_2\in[0,1]$. Then from Eq.\eqref{4.1} we have  $I(x_1,y)=D(N(x_1),y)=D(N_I(x_1),y)=D(N_I(x_2),y)=D(N(x_2),y)=I(x_2,y)$ for any $y\in[0,1]$.

(ii)$\Rightarrow$(i)  Assume that $I$ is a fuzzy implication and $N_I$ is a continuous fuzzy negation. We shall show that $I_{{D_I},{N_I}}=I$. Let $x,y\in[0,1]$. If $x\in \mbox{Ran}(\aleph_I)$, then by Eqs.\eqref{4.1}, \eqref{4.4} and \eqref{4.6} we obtain
\begin{equation*}
  I_{{D_I},{N_I}}(x,y)=D_I(N_I(x),y)=I(\aleph_I \circ N_I(x),y)=I(x,y).
\end{equation*}
If $x\notin \mbox{Ran}(\aleph_I)$, then from the continuity of $N_I$ there exists $x_0\in \mbox{Ran}(\aleph_I)$ such that $N_I(x)=N_I(x_0)$. Therefore, by Eq.\eqref{4.7} we get $I(x,y)=I(x_0,y)$ for any $y\in[0,1]$. Then
\begin{equation*}
  I_{{D_I},{N_I}}(x,y)=D_I(N_I(x),y)=D_I(N_I(x_0),y)=I(\aleph_I \circ N_I(x_0),y)=I(x_0,y)=I(x,y).
\end{equation*}
Hence we always have that $I_{{D_I},{N_I}}(x,y)=I(x,y)$ for all $x,y\in [0,1]$. Therefore, by Definition \ref{def4.2}, $I$ is a $(D_I,N_I)$-implication. On the other hand, by Eqs.\eqref{4.3} and \eqref{4.6} we have that $D_I(x,0)=I(\aleph_I(x),0)=N_I\circ\aleph_I(x)=x$ for any $x\in[0,1]$, i.e., $D_I$ is a fuzzy disjunction with a right neutral element $0$.

Finally we shall show that the representation Eq.\eqref{4.1} of $(D,N)$-implication is unique. Assume that there exist two continuous fuzzy negations $N_1$, $N_2$ and two fuzzy disjunctions $D_1$,$D_2$ which have right neutral elements 0 such that $I(x,y)=D_1(N_1(x),y)=D_2(N_2(x),y)$ for any $x,y\in[0,1]$. Then from the proof of (i)$\Rightarrow$(ii),  we have $N_1=N_I=N_2$. Now, let $x_0\in [0,1]$. Then there exists an $x_1\in[0,1]$ such that $N_I(x_1)=x_0$ since $N_I$ is a continuous fuzzy negation. Thus  $D_1(x_0,y)=D_1(N_I(x_1),y)=I(x_1,y)=D_2(N_I(x_1),y)=D_2(x_0,y)$ for any $y\in [0,1]$. Therefore, $D_1=D_2$ by the arbitrariness of $x$.
\end{proof}

Notice that in Theorem \ref{the4.1}, if $N_I$ is a strict fuzzy negation then Eq. \eqref{4.7} is naturally satisfied. Therefore, we have the following two corollaries.
\begin{corollary}\label{cor4.1}
For a function $I:[0,1]^2\rightarrow[0,1]$ the following statements are equivalent:\\
(i) $I$ is a $(D,N)$-implication in which $N$ is a strict fuzzy negation and $D$ is a fuzzy disjunction with a right neutral element $0$.\\
(ii) $I$ is a fuzzy implication and $N_I$ is a strict fuzzy negation.
\end{corollary}
\begin{corollary}
For a function $I:[0,1]^2\rightarrow[0,1]$ the following statements are equivalent:\\
(i) $I$ is a $(D,N)$-implication in which $N$ is a strong fuzzy negation and $D$ is a fuzzy disjunction with a right neutral element $0$.\\
(ii) $I$ is a fuzzy implication and $N_I$ is a strong fuzzy negation.
\end{corollary}
\begin{remark}
\emph{In Table 3 we show that the conditions in Theorem \ref{the4.1} (ii) are independent from each other, where $F$ represents the function from $[0, 1]^2$ to $[0, 1]$, i.e., $N_F(x)=F(x,0)$ for all $x\in[0,1]$.}
\begin{table}[!h]
  \centering
  \caption{The mutual independence of conditions in Theorem \ref{the4.1} (ii)}\label{3}
\begin{tabular}{|c|c|c|c|}
  \hline
  Function $F$ & $F\in\mathcal{FI}$ & Condition \eqref{4.7} & $N_F$ is continuous \\
  \hline
  $F(x,y)=I_{RS}=\begin{cases}
                   1 & \mbox{if } x\leq y, \\
                   0 & \mbox{if } x>y
                 \end{cases}$ & $\checkmark$ & $\times$ & $\times$ \\
  \hline
  $F(x,y)=\begin{cases}
             1 & \mbox{if } x=0, \\
             0.5 & \mbox{if } 0<x<0.5, \\
             y & \mbox{if } x\geq0.5
           \end{cases}$ & $\times$ & $\checkmark$ & $\times$ \\
  \hline
  $F(x,y)=\begin{cases}
            \max\{-2x+y+1,0\} & \mbox{if } x\in [0,0.5], \\
            0 & \mbox{otherwise}
          \end{cases}$ & $\times$ & $\times$ & $\checkmark$ \\
  \hline
  $F(x,y)=\begin{cases}
            1 & \mbox{if } x=0, \\
            \max\{0.5,y\} & \mbox{if } 0<x\leq0.5, \\
            y & \mbox{if } x>0.5
          \end{cases}$ & $\checkmark$ & $\checkmark$ & $\times$ \\
  \hline
  $F(x,y)=\begin{cases}
            \max\{1-2x,0\} & \mbox{if } x>y, \\
            1 & \mbox{otherwise}
          \end{cases}$ & $\checkmark$ & $\times$ & $\checkmark$ \\
  \hline
  $F(x,y)=\begin{cases}
            0 & \mbox{if } x<0.5, y>0, \\
            1 & \mbox{if } x\geq0.5, y>0, \\
            1-x & \mbox{if } y=0
          \end{cases}$ & $\times$ & $\checkmark$ & $\checkmark$ \\
  \hline
\end{tabular}
\end{table}
\end{remark}
\begin{theorem}\label{thm4.2}
 For a function $I:[0,1]^2\rightarrow[0,1]$ the following statements are equivalent:\\
(i) $I$ is a continuous $(D,N)$-implication in which $N$ is a strict fuzzy negation and $D$ is a fuzzy disjunction with a right neutral element $0$. \\
(ii) $I$ is a $(D,N_I)$-implication with a continuous fuzzy disjunction $D$ and a strict fuzzy negation $N_I$.
\end{theorem}
\begin{proof}
(i)$\Rightarrow$(ii) By Corollary \ref{cor4.1}, $N_I$ is a strict fuzzy negation. From Theorem \ref{the4.1} we know that $D_I=D$, where $D_I$ is given by Eq.\eqref{4.6}. On the other hand, by Remark \ref{rem4.3} (i) we have that $\aleph_I=N_I^{-1}$. Thus $\aleph_I$ is strict. Because $I$ is continuous and $D(x,y)=D_I(x,y)=I(\aleph_I(x),y)$ for all $x,y\in [0,1]$, we know that $D$ is continuous.

(ii)$\Rightarrow$(i) By Eq.\eqref{4.1}, $I$ is a continuous $(D,N_I)$-implication. From Corollary \ref{cor4.1},  $I$ is a $(D,N_I)$-implication in which $N_I$ is a strict fuzzy negation and $D$ is a fuzzy disjunction with a right neutral element $0$. From Theorem \ref{the4.1} we know that $N_I=N$. Therefore, $I$ is a continuous $(D,N)$-implication in which $N$ is a strict fuzzy negation and $D$ is a fuzzy disjunction with a right neutral element $0$.
\end{proof}

\begin{remark}
\emph{In Theorem \ref{thm4.2} (i), the condition that $N$ is strict cannot be moved. For instance, consider the fuzzy negation $N$ given in Remark \ref{rem4.2}. From Eq.\eqref{4.2} one can check that}
$$\aleph(x)=\begin{cases}
             -0.5x+0.5 & \emph{if } x\in(0,1], \\
             1 & \emph{if } x=0.
           \end{cases}$$
\emph{Consider the continuous fuzzy implication $I(x,y)=I_{RC}=1-x+xy$ for all $x,y\in[0, 1]$. Then by Eq.\eqref{4.6} we have  $$D(x,y)=D_I(x,y)=I(\aleph(x),y)=1-(1-y)\aleph(x)=\begin{cases}
		0.5(x+y-xy+1) & \emph{if } x\in(0,1], \\
		1-x+xy & \emph{if } x=0.
	\end{cases}$$
It is obvious that $D$ is not continuous with respect to the first variable at the point $x=0$.}
\end{remark}

As an application of Theorem \ref{thm3.2}, we finally give a characterization of a $(D,N_D)$-implication. We need the following lemma.
\begin{lemma}[\cite{QM2017}]
For a function $I:[0,1]^2\rightarrow[0,1]$ the following statements are equivalent:\\
(i) $I$ is a $(D,N)$-implication in which $N$ is a continuous (resp. strict, strong) fuzzy negation and $D$ is a commutative fuzzy disjunction with a neutral element $0$.\\
(ii) $I$ satisfies (I2) and R-CP($N_I$), and $N_I$ is a continuous (resp. strict, strong) fuzzy negation.	
\end{lemma}

\begin{theorem}
For a function $I:[0,1]^2\rightarrow[0,1]$ the following statements are equivalent:\\
(i) $I$ is a $(D,N_D)$-implication in which $N_D$ is a continuous fuzzy negation and $D$ is a commutative fuzzy disjunction with a neutral element $0$.\\
(ii) $I$ satisfies (I2) and R-CP($N_I$), and $N_I$ is a continuous (resp. strict, strong) fuzzy negation.
\end{theorem}
\section{Conclusions}
In this contribution, we mainly characterized both fuzzy negations and $(D,N)$-implications generated from fuzzy disjunctions and negations. It is worth noting that the characterization of $(D,N)$-implications shown in Ma and Zhou \cite{QM2017} is obtained through strengthening the conditions of fuzzy disjunctions, but in Theorem \ref{the4.1} (resp. Theorem \ref{thm4.2}) we just imposed a slightly strengthened condition on fuzzy implications (resp. negations).

\section*{Declaration of competing interest}
The authors declare that they have no known competing financial interests or personal relationships that could have appeared to influence the work reported in this article.

\end{document}